\documentclass[11pt]{article}
\usepackage[T1]{fontenc}
\usepackage{lmodern,microtype}
\usepackage[a4paper,margin=27mm]{geometry}
\usepackage{amsmath,amssymb,amsthm,mathtools,amscd}
\newtheorem{theorem}{Theorem}[section]
\newtheorem{proposition}[theorem]{Proposition}
\newtheorem{lemma}[theorem]{Lemma}
\theoremstyle{remark}
\newtheorem{example}[theorem]{Example}
\newtheorem{defi}{Definition}[section]
\newtheorem{remark}[theorem]{Remark}
\numberwithin{equation}{section}
\newcommand{\C}{\mathbb C}
\newcommand{\PP}{\mathbb P}

\newcommand{\OO}{\mathcal O}
\newcommand{\HH}{\mathcal H}
\newcommand{\MM}{\mathcal M}
\newcommand{\UU}{\mathcal U}
\newcommand{\II}{\mathcal Y}
\newcommand{\Levi}{\omega}

\DeclareMathOperator{\Pic}{Pic}
\DeclareMathOperator{\Hom}{Hom}
\DeclareMathOperator{\Aut}{Aut}
\DeclareMathOperator{\Lie}{Lie}
\DeclareMathOperator{\Spec}{Spec}
\DeclareMathOperator{\Tot}{Tot}
\DeclareMathOperator{\im}{im}

\DeclareMathOperator{\Hilb}{Hilb}

\title{The classification of Fano contact manifolds revisited}
\author{Baohua Fu}

\date{}
\begin{document}
\maketitle

\begin{abstract}
We present an account of the classification of Fano contact manifolds
based on the proof of OpenAI and the solution of LeBrun-Salamon conjecture by Hu, Liu and Wan. The proof
combines their constructions with earlier results on contact lines,
with some simplifications.
\end{abstract}

\section{Introduction}

Let $X$ be a smooth connected complex projective variety of
dimension $2n+1$. A {\em contact structure} on
$X$ is a corank-one vector subbundle $F\subset T_X$ for which the Lie
bracket induces an everywhere nondegenerate alternating form, called the
{\em Levi pairing},
\[
\Levi:\wedge^2F\longrightarrow L:=T_X/F,\qquad
\Levi(u,v)=[u,v]\bmod F.
\]
Here $u, v$ are local sections of $F$, and the induced form is
$\OO_X$-bilinear.
We call $L$ the {\em contact line bundle} and
write $\theta:T_X\to L$ for the quotient map. By adjunction, we have $-K_X\simeq L^{\otimes(n+1)}$. Thus $X$ is Fano if and only if $L$ is ample.

Here are classical examples of projective contact manifolds. For any smooth projective variety $Y$ of dimension at least two,
the projectivized cotangent bundle $\PP(T^*Y)$ carries a canonical
contact structure with contact line bundle $\OO_{\PP(T^*Y)}(1)$.
Another family of examples consists of adjoint varieties.
For a complex simple Lie algebra $\mathfrak g$, there is a unique closed
$G$-orbit $X_{\mathfrak g} \subset \PP(\mathfrak g)$ which is Fano contact with contact line bundle 
$\OO_{\mathbb P(\mathfrak{g})}(1)|_{X_{\mathfrak g}}$ (see Section \ref{sec:adjoint}).

By Kebekus--Peternell--Sommese--Wi\'sniewski
(\cite[Theorem~1.1]{KPSW}), together with Demailly's non-pseudoeffectivity
result (\cite[Corollaries~2--4]{Demailly}), every projective contact
manifold is either a projectivized cotangent bundle or a Fano manifold
of Picard number one. The classification of projective contact manifolds therefore reduces to the
Fano case, where the conjecture asserts that every Fano contact
manifold is an adjoint variety.

Beauville (\cite[Theorem~0.1]{Beauville}) proves this conjecture
when the map defined by $|L|$ is generically finite and the identity
component $\Aut^0(X,F)$ of the contact automorphism group is reductive.
For the three-dimensional case, see \cite{Ye} and  \cite[pp.~1864--1865]{BKK}.
Druel proves the five-dimensional case (\cite{Druel}).
Buczy\'nski--Wi\'sniewski classify Fano contact manifolds of
complex dimension at most nine when the identity component
$\Aut^0(X)$ is reductive (\cite[Theorem~1.2]{BW}).
Occhetta--Romano--Sol\'a Conde--Wi\'sniewski prove the conjecture
when $\Aut^0(X,F)$ is reductive of rank at least
$\max\{2,(n-3)/2\}$ (\cite[Theorem~6.1]{ORSW}).

There are many works on the study of Fano contact manifolds. Most importantly for this paper is 
the geometry of contact lines and conics, which has been thoroughly  studied among others by Kebekus
(\cite{Kebekus}, \cite{KebekusII}), Wi\'sniewski (\cite{Wis00}), Buczy\'nski, G. Kapustka, and M. Kapustka
(\cite{Buczynski} \cite{BKK}). 

The Riemannian counterpart is the LeBrun--Salamon conjecture
(\cite{LS}): every compact quaternionic--K\"ahler manifold with positive
scalar curvature is a Wolf space up to rescaling.
The twistor construction of Salamon and B\'erard-Bergery associates
to a compact positive quaternionic--K\"ahler $4n$-manifold a Fano
contact $(2n+1)$-manifold with a K\"ahler--Einstein metric
(\cite[Theorem~4.3 and 6.1]{Salamon}).
Conversely, a Fano contact manifold admitting a K\"ahler--Einstein
metric is a twistor space (\cite[Theorem~A]{LeBrun}).
The twistor spaces of Wolf spaces are precisely the adjoint varieties
of dimension at least three (\cite[Theorem~6.1]{Wolf}).
Thus the Fano contact classification implies
the LeBrun--Salamon conjecture.

 Poon--Salamon prove the LeBrun--Salamon
conjecture in real dimension eight (\cite{PS}), and
Buczy\'nski--Wi\'sniewski prove it in real dimensions twelve and
sixteen (\cite[Theorem~1.1]{BW}).
Brendle--Semmelmann prove it under nonnegative sectional curvature
(\cite[Corollary~2]{BS}), while the recent preprint of
Daura Serrano--Gohr--Lawn--Schedler extends the dimension bound
to $56$ (\cite[Theorem~1.2]{DGLS}).

The final classification of Fano contact manifolds is obtained in 
\cite{OpenAI}.  Hu--Liu--Wan (\cite{HLW}) prove
the LeBrun--Salamon conjecture based on the main result of Beauville (\cite{Beauville}).
We will present the proof of \cite{OpenAI} for the classification with some simplification.

\begin{theorem}[\cite{OpenAI}]\label{thm:main}
A Fano contact manifold is an adjoint variety.
\end{theorem}

Wi\'sniewski's pioneering preprint \cite{Wis00} is an important
source of inspiration for the present proof and a precursor of the
two-point approach developed in \cite{OpenAI}, where it is not cited.
 He introduced twistor conics,
considered their Hilbert component and the associated two-point
incidence over $X\times X$, and proved that its evaluation map is
surjective. The proof of \cite{OpenAI} can be viewed as a development of
Wi\'sniewski's program: it retains the geometry of degree-two
curves and two-point evaluation, but uses a stable-map
compactification and a descended section of $L\boxtimes L$ to
obtain the required global information from the diagonal by a
different argument.

We will illustrate the proof strategy in details in Section~\ref{sec:adjoint-model}. Let us compare the proof presented here with that in \cite{OpenAI}.
The proof mainly consists of four steps: (i) the construction of twistor conics. For this step, 
we use the pair of contact lines $(\ell_1, \ell_2)$ in
\cite[Propositions~3.5]{OpenAI}, but replace \cite[Lemmas~4.2-4.3,Proposition~4.4]{OpenAI} by the  vanishing of
$H^0(D,N_{D/X}^*)=0$ for  $D=\ell_1\cup\ell_2$. 
(ii) the extension of
the section $\widetilde\eta$ (which is the section $W_{21}$ in \cite{OpenAI}) across the boundary. 
For this step, we use the relative-conic construction of \cite{EF}
and the pole cancellation of
\cite[Proposition~3.3]{HLW}, which avoids normalizing the
incidence space or reducing extension to divisorial
valuations as in \cite[Section~5]{OpenAI}.
(iii) the descent of $\widetilde\eta$ to $X \times X$. For this step,  we use Hwang's finite-image arguments
\cite{Hwang} and prove the birationality of the map $q$ to produce
$\eta\in H^0(X\times X,L\boxtimes T_X)$, in place of 
\cite[Proposition~6.1, Corollary 6.2]{OpenAI}.
(iv) the proof of the homogeneity. For this,  we use directly the quadratic term
of $\eta$ along the diagnol to prove that $T_X$ is globally generated, without the scalar symmetry and cubic-term calculation used in
\cite[Proposition~7.2, Lemma~7.3 and Theorem~7.5]{OpenAI}.

Here is the organization of this paper. Section~\ref{sec:adjoint-model} illustrates the key ideas through
adjoint varieties. Section~\ref{sec:lines} recalls previous results on the geometry of contact lines
and extends the descent result of \cite{Hwang} to the contact case. Section~\ref{sec:curve} proves the existence of twistor conics.
Section~\ref{sec:section} constructs the section $\eta \in H^0(X\times X, L \boxtimes T_X)$ and then proves the homogeneity.  \vspace{0.3cm}

\begin{remark}
    Hu--Liu--Wan (\cite{HLW2}) independently give a presentation of the proof. Both arguments extend the contact contraction of the two-point tensor across nodal curves
using the pole-cancellation of \cite{HLW}. They obtain balanced
transverse maps through every point by a two-point evaluation
calculation and a propagation argument
\cite[Proposition~2.11 and Corollary~2.15]{HLW2}, whereas we
construct twistor conics by conormal vanishing.
Then they prove the irreducibility and nonlinearity of the
general VMRT and apply Hwang's theorem to descend incidence
sections whose values span each tangent space
\cite[Propositions~3.7 and~3.9]{HLW2}.
Here we adapt Hwang's argument directly to contact-line geometry,
without first establishing these VMRT properties.
The descended vector fields imply birationality of
$q:Y\to X\times X$, yielding a section
$\eta\in H^0(X\times X,L\boxtimes T_X)$ whose quadratic term
along the diagonal proves that $T_X$ is globally generated.
\end{remark}

{\em Acknowledgements:}  This paper grew out from my attempts to understand \cite{OpenAI} and \cite{HLW} through extensive  discussions with ChatGPT. 
My contribution is mainly editorial: the rewriting, reformulation and reorganization of claims and proofs; the explanation and illustration of the proof via adjoint varieties (which is the essential part of my contribution in this paper) and most importantly attributions of previous works. 

I would like to thank Zhi Jiang for organizing the conference in SCMS Shanghai from October 26 to 30, 2026, which is the driving-force for me to write down this note for the talk given there. This work was supported by the National Key Research and Development Program of China (No. 2025YFA1017302).

\section{Adjoint varieties and the proof strategy}\label{sec:adjoint-model}

\subsection{Adjoint varieties and their geometry} \label{sec:adjoint}

Let $G$ be the adjoint group of a complex simple Lie algebra
$\mathfrak g$. We choose a highest-root $\mathfrak{sl}_2$-triple
$(\mathbf{e},\mathbf{h},\mathbf{f})$ with $[\mathbf{e},\mathbf{f}]=\mathbf{h}$, $[\mathbf{h},\mathbf{e}]=2\mathbf{e}$ and $[\mathbf{h},\mathbf{f}]=-2\mathbf{f}$.
The eigenspaces of $\operatorname{ad}(\mathbf{h})$ give the contact grading
\[
\mathfrak g=\bigoplus_{i=-2}^{2}\mathfrak g_i,\qquad
\mathfrak g_{-2}=\C\mathbf{f},\quad\mathfrak g_2=\C\mathbf{e}.
\]
Let $x_0=[\mathbf{e}]$ and let $P$ be its stabilizer. The adjoint variety is the unique closed $G$-orbit in $\mathbb{P}(\mathfrak g)$, which is given by $X_{\mathfrak g}=G/P \subset \mathbb{P}(\mathfrak g)$.  Assume $\dim X_{\mathfrak g}=2n+1\geq3$. Then we have
\[
\mathfrak p=\Lie(P)=\mathfrak g_0\oplus\mathfrak g_1\oplus\mathfrak g_2
\ \text{and} \ T_{x_0}X_{\mathfrak g}\simeq\mathfrak g/\mathfrak p
\simeq\mathfrak g_{-1}\oplus\C\mathbf{f}.
\]
The $P$-invariant subspace $\mathfrak g_{-1}\subset\mathfrak g/\mathfrak p$
defines a $G$-invariant contact distribution $F_{\mathfrak g}\subset T_{X_{\mathfrak g}}$.
Its Levi pairing at $x_0$ is the Lie bracket $[u,v]=\omega_0(u,v)\mathbf{f}$, where
$\omega_0$ is nondegenerate on $\mathfrak g_{-1}$.
The contact line bundle is given by $L_{\mathfrak g}=\OO_{X_{\mathfrak g}}(1)$. Let  $\theta_{\mathfrak g}: T_{X_{\mathfrak g}}
\to L_{\mathfrak g}$ be the quotient map.

Let $\langle, \rangle$ be the scalar multiple of the Killing form normalized by  $\langle \mathbf{e},\mathbf{f}\rangle=1$.
For $c\in\mathfrak g$, we denote by $\xi_c$ the vector field on $X_{\mathfrak g}$ induced from $c$. Then for $x\in X_{\mathfrak g}$, the tangent vector $\xi_c(x) \in T_x X_{\mathfrak g}$ is given by
$\xi_c(x)=\left.\frac{d}{d\tau}\right|_{\tau=0}\exp(\tau c)\cdot x$.
For any $[b] \in X_{\mathfrak g}$, the map $\theta_{\mathfrak g}$ is given by 
\[
\theta_{\mathfrak g,[b]}(\xi_c([b]))(b)=\langle c,b\rangle,
\ \text{with} \ \theta_{\mathfrak g,x_0}(w+t\mathbf{f})(\mathbf{e})=t.
\]

We have the affine chart of $X_{\mathfrak g}$ given by the action of the opposite unipotent:
\begin{equation}\label{eq:adjoint-chart}
\Phi:\mathfrak g_{-1}\oplus\C\longrightarrow X_{\mathfrak g},\qquad
(w,t)\longmapsto[\widehat\Phi(w,t)],\qquad
\widehat\Phi(w,t)=\exp\bigl(\operatorname{ad}(w+t\mathbf{f})\bigr)\mathbf{e}.
\end{equation}
Put $Q_2(w)=\tfrac12\operatorname{ad}(w)^2\mathbf{e}\in\mathfrak g_0$.
In our notation, the parametrization of Landsberg--Manivel
\cite[\S\S2.2 and 6.1; see also \S6.3]{LM02} gives, 
\begin{equation}\label{eq:adjoint-parametrization}
\widehat\Phi(w,t)
=\mathbf{e}+[w,\mathbf{e}]-t\mathbf{h}+Q_2(w)-tw-t^2\mathbf{f}+O\bigl((w,t)^3\bigr).
\end{equation}
For types other than $A$ and $C$, the tangent directions to lines
through $x_0$ form the subvariety
$\mathcal C_{x_0}=\{[w]\in\PP(\mathfrak g_{-1})\mid Q_2(w)=0\}$
(\cite[Proposition~6.2 and Remark~6.3]{LM02}).

Adjoint varieties are key examples of Heisenberg symmetric varieties as
introduced in \cite{FH}. The higher terms in \eqref{eq:adjoint-parametrization} correspond to  contact fundamental forms of $X_{\mathfrak g} \subset \mathbb{P}(\mathfrak g)$ (see \cite[Lemma~6.6 and Remark~6.7]{FH}).

\begin{remark}
For $\mathfrak{g}$ not of type $A$ or $C$, 
Mok's recognition theorem (\cite[Main Theorem]{MokVMRT}) proves that
a Fano manifold of Picard number one is isomorphic to an
adjoint variety if its VMRT at a general point is
projectively equivalent to $\mathcal C_{x_0} \subset \PP(\mathfrak g_{-1})$. Thus a natural approach to prove Theorem \ref{thm:main} is to determine the VMRT of a Fano contact manifold.
Hwang emphasized this approach (\cite[Section~1, Theorem~1.4 and the following paragraph]{HwangCone}).
His work on Legendrian cone structures
(\cite[Theorem~1.6]{HwangCone}) and later criteria for
higher-dimensional subadjoint Legendrian varieties
(\cite[Theorems~1.2--1.3]{HwangSubadjoint}) are recent advances of this approach.
\end{remark}
\subsection{A section from two points on a rational curve}

 Consider  a Fano contact manifold $(X,F,L,\theta)$ of dimension $2n+1$.
Let $h:\Sigma\simeq\PP^1\to X$ be a nonconstant morphism.
For $a\in\Sigma$, evaluation at $a$ of the trivial line bundle
$T_\Sigma(-2a)$ gives a canonical isomorphism
\[
\rho_a:H^0(\Sigma,T_\Sigma(-2a))\xrightarrow{\sim}
T_a\Sigma\otimes\operatorname{Sym}^2T_a^*\Sigma
\simeq T_a^*\Sigma.
\]
In local coordinate $u$ of $\Sigma$ with $u(a)=0$, $\rho_a$ sends
$c u^2\partial_u$ to $c\,du|_a$. It follows that $\rho_a^{-1}$ identifies a cotangent vector at $a$ to a vector field on $\Sigma$ which vanishes at $a$ to order two.

For $b\in\Sigma$, let $x=h(a)$ and $y=h(b)$. Consider the following composition
\[
T_a^*\Sigma\xrightarrow{\rho_a^{-1}}
H^0(\Sigma,T_\Sigma(-2a))\xrightarrow{\mathrm{ev}_b}
T_b\Sigma\xrightarrow{dh_b}T_yX.
\]
Tensoring with $L_x$, we get a map $L_x\otimes T_a^* \Sigma \to L_x \otimes T_yX$.
The map $\theta_x\circ dh_a: T_a\Sigma \to L_x$ corresponds to a vector in $L_x\otimes T_a^* \Sigma$. Its image under the previous map gives an element in $L_x\otimes T_yX$, denoted by
$\widetilde\eta([h;a,b])$.

Let us denote by $Y_0$ the set of all 
triples $[h;a,b]$ with $h: \Sigma \to X$ non-constant and $a, b \in \Sigma$. Then we have  a natural
two-point evaluation map

$$q: Y_0 \to X \times X, \qquad [h;a,b] \mapsto (h(a),h(b)). $$ The preceding construction gives
a section $\widetilde\eta \in H^0(Y_0, q^*(L\boxtimes T_X))$, which is the key construction in \cite{OpenAI} (denoted by $W_{21}$ there).

A first question is when $\widetilde\eta$ could be nonzero. 
If $d=\deg h^*L$, then the contact differential $\theta \circ dh$ satisfies
$$
\theta\circ dh\in H^0(\Sigma,\Omega_\Sigma^1\otimes h^*L)
\simeq H^0(\PP^1,\OO(d-2)).
$$
If $d=1$ (namely a contact line), then we always get zero. 
Thus the first interesting case is $d=2$ and in this case, the contact differential
$\theta\circ dh$ is either zero or an isomorphism.
Moreover, to descend $\widetilde\eta$ to $X\times X$, we need the
two-point evaluation $q$ to be dominant. This leads us naturally to consider twistor conics
(Wi\'sniewski \cite[p.~1]{Wis00}).

\begin{defi}
   A smooth rational
curve $C\subset X$ is called a {\em twistor conic} if
\[
L|_C\simeq\OO_C(2),\qquad
N_{C/X}\simeq\OO_C(1)^{\oplus2n} \quad \text{and} \quad
\theta|_{T_C}:T_C\xrightarrow{\sim}L|_C.
\]
\end{defi}

Among smooth rational curves of $L$-degree two, these conditions
define an open subset of the Hilbert scheme. 
Let us consider  the simplest example:
$X=\Sigma=\PP^1$ and $\theta:T_\Sigma\to L$ is an isomorphism.
Then the variety $Y_0$ is nothing else but just $\Sigma \times \Sigma$. We will describe the section $\tilde{\eta}$ just constructed.
Let $\Delta_\Sigma\subset\Sigma\times\Sigma$
be the diagonal. Consider the line bundle
$\mathcal L_\Sigma=(T_\Sigma\boxtimes T_\Sigma)(-2\Delta_\Sigma)$.  The identification $N_{\Delta_\Sigma/(\Sigma\times\Sigma)}\simeq T_\Sigma$ (induced by $(v_a,v_b)\mapsto v_b-v_a$) gives a natural trivialization
\[
\mathcal L_\Sigma|_{\Delta_\Sigma}
\simeq T_\Sigma^{\otimes2}\otimes(\Omega_\Sigma^1)^{\otimes2}
\simeq\OO_\Sigma.
\]
Restriction to the diagonal identifies
$H^0(\Sigma\times\Sigma,\mathcal L_\Sigma)$ with $\C$.
Define $\psi_\Sigma$ to be the image in
$H^0(\Sigma\times\Sigma,T_\Sigma\boxtimes T_\Sigma)$ of the
unique section whose restriction is $1$ under this identification. This is exactly
the section $\widetilde\eta \in H^0(\Sigma \times \Sigma, L \boxtimes T_X)$.

For an affine coordinate $u$ induced by an isomorphism
$\Sigma\simeq\PP^1$, let $u_a,u_b$ be the pullbacks of $u$ to the two factors.
Then this definition gives the following two-point tensor (which is denoted by $q$ in \cite[Section~3.1, (3.1)]{HLW} and by $b$ in
\cite[Section~5.2, (5.5)]{OpenAI}).

\begin{equation}\label{eq:two-point-tensor}
\psi_\Sigma=(u_b-u_a)^2\partial_{u_a}\otimes\partial_{u_b},
\end{equation}

For the general case, the explicit formula for $\widetilde\eta$ is
\begin{equation}\label{eq:curvewise-section}
\begin{aligned}
\widetilde\eta([h;a,b])
&=\bigl((\theta_{h(a)}\circ dh_a)\otimes dh_b\bigr)
(\psi_\Sigma(a,b))\\
&=(u(b)-u(a))^2\theta_{h(a)}\bigl(dh_a(\partial_u)\bigr)
\otimes dh_b(\partial_u).
\end{aligned}
\end{equation}
The construction is invariant under reparametrization and works in
families.

\subsection{The global section on an adjoint variety}

Let us return to the adjoint variety $X_{\mathfrak g}$.
The tautological inclusion $L_{\mathfrak g}^{-1} \subset \mathfrak{g} \otimes \mathcal{O}_{X_{\mathfrak g}}$ and the infinitesimal action $\mathfrak g \otimes \mathcal{O}_{X_{\mathfrak g}} \to T_{X_{\mathfrak g}}$ give
\[
{\rm pr}_1^*L_{\mathfrak g}^{-1} \hookrightarrow 
\mathfrak g\otimes\OO_{X_{\mathfrak g}\times X_{\mathfrak g}} \longrightarrow
{\rm pr}_2^*T_{X_{\mathfrak g}}
\]
The composition defines the $G$-invariant regular section (up to a sign)
\begin{equation}\label{eq:adjoint-tensor}
\eta_{\mathfrak g}\in H^0(X_{\mathfrak g}\times X_{\mathfrak g},L_{\mathfrak g}\boxtimes T_{X_{\mathfrak g}}),\qquad
\eta_{\mathfrak g}([c],x)(c)=-\xi_c(x).
\end{equation}

The adjoint variety $X_{\mathfrak g}$ contains in fact many twistor conics (which is already known from \cite{Salamon}, see also \cite{Kwon} for a systematic study of conics on $X_{\mathfrak g}$). 
We are going to show that $\eta_{\mathfrak g}$ agrees with $\tilde{\eta}$ on the $G$-translates of  twistor conics.
In our notation, the highest-root triple $(\mathbf{e}, \mathbf{h}, \mathbf{f})$ gives an $\mathrm{SL}_2$-equivariant embedding
\begin{equation}\label{eq:adjoint-conic}
j:\PP^1\longrightarrow X_{\mathfrak g},\qquad
[s:t]\longmapsto[s^2\mathbf{e}-st\mathbf{h}-t^2\mathbf{f}].
\end{equation}
Write $C=j(\PP^1)$ and $u=t/s$, so $j(0)=x_0$ and
$j(u)=[\operatorname{Ad}(\exp(u\mathbf{f}))\mathbf{e}]$.
Since $j^*L_{\mathfrak g}\simeq\OO(2)$ and
$\theta_{\mathfrak g,x_0}(dj_0(\partial_u))(\mathbf{e})=1$, the contact differential $\theta_{\mathfrak g}(dj)$ is
an isomorphism. 
%The root Borel acts on the normal fibre
%$T_{x_0}X_{\mathfrak g}/T_{x_0}C\simeq\mathfrak g_{-1}$ with weight $-1$ and trivial unipotent action, giving 
Moreover, $N_{C/X_{\mathfrak g}}\simeq\OO_C(1)^{\oplus2n}$.
Thus $C$ is a twistor conic.

Under this $\mathrm{SL}_2$-action, $\exp(\tau\mathbf{e})$ sends
 $u$ to $u/(1+\tau u)$.
It follows that
$-\xi_{\mathbf{e}}(j(u))=dj_u(u^2\partial_u)$.
Together with $\theta_{\mathfrak g,x_0}(dj_0(\partial_u))(\mathbf{e})=1$, this proves
that $\eta_{\mathfrak g}$ agrees with \eqref{eq:curvewise-section} when $a=0$.
By the  $\mathrm{SL}_2$-equivariance, we have the identity for every
$a,b\in\PP^1$,
\begin{equation}\label{eq:adjoint-restriction}
\eta_{\mathfrak g}(j(a),j(b))
=(\theta_{\mathfrak g,j(a)}\circ dj_a\otimes dj_b)
\bigl(\psi_{\PP^1}(a,b)\bigr).
\end{equation}

By $G$-equivariance, \eqref{eq:adjoint-restriction} holds on
every translate of $C$. For types other than $A$ or $C$, every twistor conic is such a translate
by \cite[Lemma~4.2(3)]{Kwon}. Thus $q^*\eta_{\mathfrak g}=\widetilde\eta$ on the locus of
two-pointed twistor conics in these cases.

\subsection{The quadratic term and homogeneity}

We are now to recover the homogeneity of $X_{\mathfrak g}$ from
$\eta_{\mathfrak g}$. Denote by
$\mathbf{s}_{x_0}=\eta_{\mathfrak g}(-,x_0)$ the restriction of $\eta_{\mathfrak g}$ to
$X_{\mathfrak g}\times\{x_0\}$.
In the chart \eqref{eq:adjoint-chart}, we trivialize
$\Phi^*L_{\mathfrak g}$ by evaluation on $\widehat\Phi(w,t)$.
Recall that $T_{x_0}X_{\mathfrak g}\simeq\mathfrak g/\mathfrak p \simeq \mathfrak g_{-1}\oplus\mathfrak g_{-2}$. Let
${\rm pr}_{-}$ denote the projection onto this
summand along $\mathfrak p$. By \eqref{eq:adjoint-tensor},
\[
\mathbf{s}_{x_0}(\Phi(w,t))\bigl(\widehat\Phi(w,t)\bigr)
=-{\rm pr}_{-}\bigl(\widehat\Phi(w,t)\bigr).
\]
In \eqref{eq:adjoint-parametrization}, the term
$\mathbf{e} + [w,\mathbf{e}]-t\mathbf{h} + Q_2(w)$
belongs to  $\mathfrak p$, hence we have
\[
-{\rm pr}_{-}\bigl(\widehat\Phi(w,t)\bigr)
=tw+t^2\mathbf{f}+O\bigl((w,t)^3\bigr).
\]

This implies that $\mathbf{s}_{x_0}$ vanishes to order at least two at $x_0$.
Via $d\Phi_{(0,0)}$, the pair $(w,t)$ corresponds to
$v=w+t\mathbf{f}\in T_{x_0}X_{\mathfrak g}$. Since
$\theta_{\mathfrak g,x_0}(v)(\mathbf{e})=t$, the quadratic
term $t(w+t\mathbf{f})$ in the chosen trivialization
corresponds to $\theta_{\mathfrak g,x_0}(v)\otimes v$.
By equivariance, for every $x\in X_{\mathfrak g}$
the section $\mathbf{s}_x=\eta_{\mathfrak g}(-,x)$ vanishes to order
at least two at $x$ with quadratic term
$v\mapsto\theta_{\mathfrak g,x}(v)\otimes v$.

These vanishings give order at least two along the diagonal
$\Delta_{\mathfrak g}$. 
We compute its quadratic term $[\eta_{\mathfrak g}]_{2,x} \in {\rm Sym}^2 T_x^*X_{\mathfrak g} \otimes L_{\mathfrak g,x} \otimes T_xX_{\mathfrak g}$.
Under
$N_{\Delta_{\mathfrak g}/(X_{\mathfrak g}\times X_{\mathfrak g})}\simeq T_{X_{\mathfrak g}}$ (defined by
$(v_1,v_2)\mapsto v_2-v_1$),  the direction $v$ in the first slice
represents $-v$. Since
$\theta_{\mathfrak g,x}(-v)\otimes(-v)=\theta_{\mathfrak g,x}(v)\otimes v$, we obtain
\begin{equation}\label{eq:adjoint-quadratic}
\begin{gathered}
\eta_{\mathfrak g}\in H^0\!\left(X_{\mathfrak g}\times X_{\mathfrak g},
\mathcal I_{\Delta_{\mathfrak g}}^2\otimes(L_{\mathfrak g}\boxtimes T_{X_{\mathfrak g}})\right),\\
[\eta_{\mathfrak g}]_{2,x}(v)=\theta_{\mathfrak g,x}(v)\otimes v.
\end{gathered}
\end{equation}

Let us show that $X_{\mathfrak g}$ is homogeneous using this section $\eta_{\mathfrak g}$.
Let $E_x=\im(\operatorname{ev}_x:H^0(X_{\mathfrak g},T_{X_{\mathfrak g}})\to T_xX_{\mathfrak g})$.
By construction, $\mathbf{s}_x=\eta_{\mathfrak g}(-,x)$ takes values in
$L_{\mathfrak g}\otimes E_x$. Its quadratic term therefore gives
\[
\theta_{\mathfrak g,x}(v)\otimes v=[\mathbf{s}_x]_{2,x}(v)\in L_{\mathfrak g,x}\otimes E_x
\qquad(v\in T_xX_{\mathfrak g}).
\]
For $v\notin\ker\theta_{\mathfrak g,x}$, contraction by a suitable element
of $L_{\mathfrak g,x}^*$ gives $v\in E_x$. Such vectors span $T_xX_{\mathfrak g}$, so
$E_x=T_xX_{\mathfrak g}$. It follows that $T_{X_{\mathfrak g}}$ is globally generated, so 
every $G$-orbit is open and $X_{\mathfrak g}$ is homogeneous. 

\subsection{The proof strategy}
For a Fano contact manifold $(X,L)$ with $\Pic(X)=\mathbb Z[L]$,
it is now clear what remains to be done:  (i) existence of  twistor conics; (ii)  compactification of the variety $Y_0$ of all triples $[h; a,b]$; (iii) extension of $\widetilde\eta$ to the compactification; (iv) descent of $\widetilde\eta$ to $\eta \in H^0(X\times X, L \boxtimes T_X)$ and deduction of the homogeneity from the quadratic terms of $\eta$ along the diagonal. 

For (i), we will pick two contact lines through a general point in  special position and then deform their union to a smooth conic. For (ii), a  natural choice to  compactify $Y_0$ is to take stable maps  of rational curves with two marked points. Then the boundary would contain chains of rational curves with contracted components, which would create more difficulties to extend $\tilde{\eta}$ across the boundary. Instead, the idea of \cite{OpenAI} is to use the fiber product of the universal family of the moduli of 
stable maps of rational curves without marked points. Then the boundary will only consist of stable maps with two contact lines and then the local computation can be carried out to prove the extension. For the descent, we will   use the descent arguments from \cite{Hwang} and also prove that the morphism $q$ is in fact birational (inspired from the adjoint case of \cite{Kwon}).

\section{Contact lines and descent of sections}\label{sec:lines}\label{sec:descent}

For the Fano contact manifold $(X, L)$, \cite[Proposition~2.3]{BKK} implies that either $\Pic(X)=\mathbb Z[L]$, or
$(X,L)\simeq\bigl(\PP^{2n+1},\OO(2)\bigr) \text{or}
\bigl(\PP(T^*\PP^{n+1}),\OO(1,1)\bigr).$
The latter two varieties are the adjoint varieties of types $C_{n+1}$ and $A_{n+1}$.
%On projective space, a contact structure is determined by a
%nondegenerate alternating form on $\C^{2n+2}$, unique up to a linear
%change of coordinates. On the incidence variety, the contact structure
%is unique \cite[Proposition~2.3]{BKK}.
For the rest of the paper, we assume $\Pic(X)=\mathbb Z[L]$. 
%Then $n\geq2$ by \cite[Theorem~2.5]{BKK}.

A {\em contact line} is an integral rational curve $\ell\subset X$
with $L\cdot\ell=1$. Such curves exist by
\cite[Section~2.3]{Kebekus}. 
Fix an irreducible component $\HH$ of the Chow space of contact
lines, and put $\HH_x=\{[\ell]\in\HH:x\in\ell\}$.
The space $\HH$ is projective and it covers $X$ by
\cite[Corollary~5.9]{BKK}. Here are some known facts on the geometry of contact lines.

\begin{proposition}\label{prop:standard-lines}
For a general point $x\in X$, the following hold.
\begin{enumerate}
\renewcommand{\theenumi}{\roman{enumi}}
\renewcommand{\labelenumi}{(\theenumi)}
\item\label{item:standard-splitting}
Every line $\ell$ with $[\ell]\in\HH_x$ is smooth. Two distinct
such lines meet only at $x$ and have distinct tangent directions
there. Moreover, if we write $\OO_\ell(1)=L|_\ell$, then we have
\begin{equation}\label{eq:standard}
T_X|_\ell\simeq\OO_\ell(2)\oplus
\OO_\ell(1)^{\oplus(n-1)}\oplus\OO_\ell^{\oplus(n+1)}.
\end{equation}

\item\label{item:positive-subbundle}
The bundles
\[
P_\ell=\im\!\left(
H^0(\ell,T_X|_\ell\otimes\OO_\ell(-1))\otimes\OO_\ell(1)
\longrightarrow T_X|_\ell\right),\qquad
Q_\ell=T_X|_\ell/P_\ell
\]
satisfy $P_\ell\subset F|_\ell$ and
$Q_\ell\simeq\OO_\ell^{\oplus(n+1)}$.
The subbundle $P_\ell$ is Lagrangian for $\Levi$.

\item\label{item:vmrt}
Let $\mathcal C_x\subset\PP(F_x)$ be the reduced variety of
tangent directions to members of $\HH_x$, and let
$\widehat{\mathcal C}_x\subset F_x$ be its affine cone.
Then $\mathcal C_x$ is smooth, possibly disconnected, Legendrian
of dimension $n-1$, and spans $\PP(F_x)$.
For a general $[\ell]\in\HH_x$ and every $0\ne v\in T_x\ell$, we have
%\begin{equation}\label{eq:tangent-cone}
$T_v\widehat{\mathcal C}_x=(P_\ell)_x.$
%\end{equation}

\item\label{item:vmrt-secants}
If $\mathcal C_x$ is irreducible, its secant variety is
$\PP(F_x)$. If $\mathcal C_1,\mathcal C_2$ are distinct irreducible components of
$\mathcal C_x$, their join $J(\mathcal C_1,\mathcal C_2)$ is $\PP(F_x)$.
Here $J(\mathcal C_1,\mathcal C_2)$ is the closure of the union of lines joining
a point of $\mathcal C_1$ to a distinct point of $\mathcal C_2$.
\end{enumerate}
\end{proposition}

\begin{proof}
The smoothness and splitting in \textup{(i)} follow from
\cite[Proposition~3.3 and Lemma~3.5]{Kebekus}. The assertions about
tangent directions and intersections follow from
\cite[Theorems~4.1 and~5.1]{KebekusII}.  The intersection
assertion is \cite[Theorem~3.1(iii)]{Buczynski}.
The inclusion and Lagrangian assertion in \textup{(ii)} are
\cite[Lemma~4]{Wis00}; the description of $Q_\ell$ follows from
\eqref{eq:standard}.

In \textup{(iii)}, smoothness, the Legendrian property and dimension
are given by \cite[Theorem~3.1(ii)--(iii)]{Buczynski},
and spanning is \cite[Theorem~4.4]{Kebekus}.
The tangent-space description  follows from 
\cite[Proposition~2.1]{HR}.

If $\mathcal C_x$ is irreducible, its secant variety fills $\PP(F_x)$
by \cite[proof of Proposition~17(2)]{LM}, since it spans $\PP(F_x)$.
Otherwise its irreducible components are disjoint by smoothness,
and \cite[Proposition~3.5]{Buczynski} gives the join assertion
in \textup{(iv)}.
\end{proof}

A vector field $\xi$ is {\em contact} if $[\xi,F]\subset F$.
By \cite[Proposition~2.1]{LeBrun}, every local section $s$ of $L$
determines a unique contact vector field $\kappa_s$ satisfying
$\theta(\kappa_s)=s$.  By the uniqueness of the contact structure \cite[Corollary~4.5]{Kebekus}, we have $H^0(X,F)=0$ and every global vector field on $X$ is contact, which can be generalized as follows:

\begin{lemma}\label{lem:F-vanishing}
Let $q:Y\to X$ be a surjective morphism from an integral
projective variety. Then $H^0(Y,q^*F)=0$. 
\end{lemma}

\begin{proof}
We can use Stein factorization and the projection formula to reduce the assertion
to the case $q$ is finite. Let $s \in H^0(Y,q^*F)$ and choose a general point
$y\in Y$ such that $q$ is \'etale at $y$ and
Proposition~\ref{prop:standard-lines} holds at $x=q(y)$.
For every $[\ell]\in\HH_x$, an irreducible component
of $q^{-1}(\ell)$ through $y$ dominates $\ell$. Its normalization
gives a finite morphism $g:\widetilde\ell\to\ell$.

The Levi pairing with $T_\ell\subset F|_\ell$ induces
$F|_\ell\to\Omega^1_\ell\otimes L|_\ell\simeq\OO_\ell(-1)$.
The image of $s|_{\widetilde\ell}$ is zero because
$\deg g^*\OO_\ell(-1)<0$. At a point over $y$, this gives
$\Levi_x(s(y),T_x\ell)=0$.
The spaces $T_x\ell$ span $F_x$ by Proposition~\ref{prop:standard-lines} and $\Levi_x$ is nondegenerate,
so $s(y)=0$, which gives $s=0$.
\end{proof}

The following is a contact counterpart of Hwang's result
\cite[Theorem~1.3]{Hwang}. His theorem assumes that the variety of
minimal rational tangents is finite or has a nonlinear irreducible
component, which is not immediate for Fano contact manifolds.
Our proof is an adaption of his arguments
\cite[Propositions~2.3--2.4]{Hwang}.

\begin{lemma}\label{lem:descent}
Let $q:Y\to X$ be a surjective morphism from an integral
projective variety. Then the pull-back induces an isomorphism
\[
q^*:H^0(X,T_X)\xrightarrow{\sim}H^0(Y,q^*T_X).
\]
\end{lemma}

\begin{proof}
As $q$ is surjective, $q^*$ is injective. It remains to prove the surjectivity.
Let $s\in H^0(Y,q^*T_X)$ and put $\chi=(q^*\theta)(s) \in H^0(Y, q^*L)$.
Let $T\subset\Tot(L)$ be the reduced image of
$y\mapsto(q(y),\chi(y))$.
The image $T$ is proper because $Y$ is projective.
Since $T\to X$ is also affine, it is finite and $T$ is projective.
If its degree is greater than one, the finite-cover lemma of
Hwang--Mok \cite[Lemma~4.2]{HM} (see also 
\cite[Proposition~2.2]{Hwang}) gives an integral component
$\Gamma\subset T|_\ell$ of degree greater than one over a general
member $\ell$ of $\HH$. Choose an integral component
$D_\ell\subset Y\times_T\Gamma$ dominating $\Gamma$, with
projections $r_\ell:D_\ell\to Y$ and $p_\ell:D_\ell\to\ell$.
The variety $D_\ell$ is projective, and $r_\ell^*s$ is a section
of $p_\ell^*(T_X|_\ell)$.

Since $H^0(D_\ell,\OO_{D_\ell})=\C$, the image of $r_\ell^*s$
in the trivial bundle $p_\ell^*Q_\ell$ is constant.
Since $\theta|_\ell$ factors through $Q_\ell$, there is
$t_\ell\in H^0(\ell,L|_\ell)$ with
$r_\ell^*\chi=p_\ell^*t_\ell$.
The image of $D_\ell$ in $\Tot(L|_\ell)$ is therefore the graph
of $t_\ell$. It is also $\Gamma$, contradicting
$\deg(\Gamma/\ell)>1$.

Thus $T\to X$ is finite birational, hence an isomorphism. Its inclusion in $\Tot(L)$ defines
$\bar s\in H^0(X,L)$ with $\chi=q^*\bar s$.
By Lemma~\ref{lem:F-vanishing},
$s-q^*\kappa_{\bar s}\in H^0(Y,q^*F)=0,$
which proves the surjectivity.
\end{proof}

\section{Twistor conics}\label{sec:curve}

Let $H\subset\Hilb(X)$ be the open subscheme parametrizing smooth
rational curves of $L$-degree two, with universal family
$p_H:U_H\to H$ and evaluation $e_H:U_H\to X$.
For $[C]\in H$, the normal bundle has rank and degree $2n$, hence $N_{C/X}\simeq\OO_C(1)^{\oplus2n}$ if and only if
$H^0(C,N_{C/X}^*)=0$. This condition and the condition that
$\theta|_{T_C}$ being an isomorphism define the open subscheme
$H^{\mathrm{tw}}\subset H$ of twistor conics.

For $x\in X$, the projection $p_H$ identifies $e_H^{-1}(x)$
with the closed subscheme $H_x\subset H$ of curves containing $x$.
We have the tangent morphism
\[
\tau_x:H_x\longrightarrow\PP(T_xX),\qquad
[C]\longmapsto[T_xC].
\]

\begin{lemma}\label{lem:twistor-tangent-map}
Assume that $[C]\in H_x$ satisfies
$N_{C/X}\simeq\OO_C(1)^{\oplus 2n}$.
Then $H_x$ is smooth of dimension $2n$ at $[C]$ and
$\tau_x$ is \'etale at $[C]$.
\end{lemma}

\begin{proof}
The first part is well known (see for example \cite[Section~3, p.1065]{IN}). Moreover $T_{[C]} H_x \simeq H^0(C,N_{C/X}(-x))$
and the tangent map is given by the evaluation:
\[
d_{[C]}\tau_x:H^0(C,N_{C/X}(-x))\longrightarrow N_{C/X, x} \otimes T_x^*C
\simeq T_{[T_xC]}\PP(T_xX).
\]
Its kernal is  $\ker(d_{[C]}\tau_x)=H^0(C,N_{C/X}(-2x))=0$, which proves that $\tau_x$ is 
\'etale at $[C]$ as both varieties are smooth of dimension $2n$.
\end{proof}

\begin{lemma}\label{lem:conormal-vanishing}
Let $x\in X$ be a general point and let
$[\ell_1],[\ell_2]\in\HH_x$ such that
$(P_{\ell_1})_x\oplus(P_{\ell_2})_x=F_x$ and $\Levi_x(T_x\ell_1,T_x\ell_2)\ne0$.
Then $D:=\ell_1\cup\ell_2$ satisfies $H^0(D,N_{D/X}^*)=0$.
\end{lemma}

\begin{proof}

By Proposition~\ref{prop:standard-lines}(\ref{item:standard-splitting}), 
$D$ is nodal and a local complete intersection, hence $N_{D/X}^*$ is locally free and
the conormal map $N_{D/X}^*\to\Omega_X^1|_D$ is injective.
As $H^0(\ell_i,P_{\ell_i}^*)=0$,  we have
$H^0(\ell_i,\Omega_X^1|_{\ell_i})=H^0(\ell_i,Q_{\ell_i}^*)$. Since $Q_{\ell_i}$ is trivial,
evaluation at $x$ identifies this space with
$\operatorname{Ann}((P_{\ell_i})_x)$. The two restrictions of a
section on $D$ agree at $x$, so evaluation gives
\[
H^0(D,\Omega_X^1|_D)\simeq
\operatorname{Ann}((P_{\ell_1})_x)\cap
\operatorname{Ann}((P_{\ell_2})_x)=\operatorname{Ann}(F_x).
\]

Suppose that $0\ne\zeta\in H^0(D,N_{D/X}^*) \subset H^0(D,\Omega_X^1|_D)$. By the preceding
isomorphism, $0\ne\zeta_x\in\operatorname{Ann}(F_x)$. Let $\mathcal I_D$ be the
ideal sheaf of $D$. Since $N_{D/X}^*=\mathcal I_D/\mathcal I_D^2$,
on an affine neighbourhood $U$ of $x$ we can represent $\zeta$
by a function $\phi\in\mathcal I_D(U)$. Thus
$d\phi|_{D\cap U}=\zeta|_{D\cap U}$. Choose a local frame of $L$
so that $\theta$ is represented on $U$ by a one-form $\alpha$
with $\alpha_x=d\phi_x$, and put $\beta=\alpha-d\phi$.
Then $\beta_x=0$ and $d\beta=d\alpha$.

Choose a nonzero $v_i\in T_x\ell_i$. As 
$T_X|_{\ell_i}$ is globally generated, there exists  $w_i \in H^0(\ell_i, T_X|_{\ell_i})$ with
$w_i(x)=v_{3-i}$. 
Recall that $H^0(\ell_i, T_X|_{\ell_i})$ parametrises first-order deformations of $\ell_i \subset X$,
hence $w_i$ defines a first-order deformation $f_i:\ell_i\times S\to X$, where
$S=\Spec(\mathbb C[\epsilon]/(\epsilon^2))$. Denote by
$p_i:\ell_i\times S\to\ell_i$  the projection. The relative
differential $d_Sf_i:p_i^*T_{\ell_i}\to f_i^*T_X$ gives a section
\[
(f_i^*\theta)\circ d_Sf_i
\in H^0(\ell_i\times S,p_i^*\Omega^1_{\ell_i}\otimes f_i^*L).
\]
Since $H^1(\ell_i,\mathcal O_{\ell_i})=0$, the line bundle
$f_i^*L$ is isomorphic to $p_i^*(L|_{\ell_i})$. Thus
$p_i^*\Omega^1_{\ell_i}\otimes f_i^*L \simeq p_i^*\mathcal O_{\ell_i}(-1)$ and
$H^0(\ell_i\times S,p_i^*\Omega^1_{\ell_i}\otimes f_i^*L)=0$. Hence
$(f_i^*\theta)\circ d_Sf_i=0$.

By the definition of the first-order deformation associated with
$w_i$, for every local function $g$ we have
$f_i^*g=g|_{\ell_i}+\epsilon\,dg(w_i)$.
Since $\phi|_{\ell_i}=0$ and $d\phi|_{\ell_i}=\zeta|_{\ell_i}$,
this gives $f_i^*\phi=\epsilon\,\zeta(w_i)$ near $x$.
The regular function $\zeta(w_i)$ on $\ell_i\simeq\PP^1$ is constant, which
 vanishes at $x$ because $w_i(x)=v_{3-i}\in F_x$ and $\zeta_x$
annihilates $F_x$. Thus $f_i^*\phi=0$. Differentiating along
$\ell_i$ gives $(f_i^*d\phi)\circ d_S f_i=0$. In the chosen frame
of $L$ we also have $(f_i^*\alpha)\circ d_S f_i= (f_i^*\theta)\circ d_S f_i=0$. Hence
$(f_i^*\beta)\circ d_S f_i=0$.

Along $x\times S$, the deformation $f_i$ moves $x$ in the
direction $w_i(x)=v_{3-i}$, while $d_Sf_i(v_i)$ reduces to $v_i$
modulo $\epsilon$. Write $(\partial_v\beta)_x$ for the
directional derivative of the coefficients of $\beta$ in local
coordinates at $x$. Since $\beta_x=0$, the coefficient of
$\epsilon$ in
$(f_i^*\beta)\circ d_Sf_i(v_i)$ is
$(\partial_{v_{3-i}}\beta)_x(v_i)$.
The pullback form vanishes, so
$(\partial_{v_{3-i}}\beta)_x(v_i)=0$ for $i=1, 2$.
It follows that
\[
(d\alpha)_x(v_1,v_2)
=(d\beta)_x(v_1,v_2)
=(\partial_{v_1}\beta)_x(v_2)
 -(\partial_{v_2}\beta)_x(v_1)=0,
\]
 which  contradicts
$\Levi_x(v_1,v_2)\ne0$ since the restriction of
$(d\alpha)_x$ to $F_x$ is $-\Levi_x$.

\end{proof}

\begin{proposition}\label{prop:curve}
The variety $X$ contains a twistor conic.
\end{proposition}

\begin{proof}

Take a general point $x\in X$. If $\mathcal C_x$ is irreducible,
put $A_1=A_2=\mathcal C_x$. Otherwise take two distinct
irreducible components $A_1,A_2$ of $\mathcal C_x$.
By Proposition~\ref{prop:standard-lines}(\ref{item:vmrt-secants}),
addition of the affine cones $\widehat A_1\times\widehat A_2\to F_x$
is dominant. At a general pair of nonzero vectors
$v_i\in T_x\ell_i$ with $[v_i]\in A_i$, its differential identifies
$(P_{\ell_1})_x\oplus(P_{\ell_2})_x$ with $F_x$
by Proposition~\ref{prop:standard-lines}(\ref{item:vmrt}).
If $\Levi_x(v_1,v_2)$ vanished identically on the product of cones,
their tangent spaces would be orthogonal. This is impossible for
complementary Lagrangian subspaces. Thus we may choose the pair with
\begin{equation}\label{eq:pair}
(P_{\ell_1})_x\oplus(P_{\ell_2})_x=F_x,\qquad
\Levi_x(v_1,v_2)\ne0.
\end{equation}
For $D=\ell_1\cup\ell_2$,
Lemma~\ref{lem:conormal-vanishing} gives $H^0(D,N_{D/X}^*)=0$.

A section of $(N_{D/X}^*|_{\ell_i})(-x)$ vanishes at $x$,
so it glues to the zero section on $\ell_{3-i}$ to get a section of $H^0(D,N_{D/X}^*)=0$.
Hence we have
$H^0\bigl(\ell_i,(N_{D/X}^*|_{\ell_i})(-x)\bigr)=0.$
It follows that every summand
of $N_{D/X}|_{\ell_i}$ has nonnegative degree. Thus $N_{D/X}|_{\ell_i}$ is globally generated with vanishing $H^1$.
Their sections can be chosen with any prescribed common value
at $x$, so the normalization sequence for $D$ gives
\[
H^1(D,N_{D/X})=0 \quad \text{and} \quad 
H^0(D,N_{D/X})\longrightarrow (N_{D/X})_x
\quad\text{surjective}.
\]
It follows that the Hilbert scheme ${\rm Hilb}(X)$ is smooth at $[D]$. 
In local coordinate the ideal of $D$ is $(uv, z_1, \cdots, z_{2n-1})$, so the $uv$-coordinate in 
$N_{D/X, x}$ controls first-order smoothing of the node. The surjection above therefore provides a global tangent vector that smooths the node.
 Choose a smooth
curve $B$ through $[D]$ in ${\rm Hilb}(X)$ with this tangent
direction, and write $\pi_B:\mathcal D_B\to B$ for the family
and $\mu:\mathcal D_B\to X$ for evaluation.

Let $p\in\mathcal D_B$ be the node of the central fibre, so that
$\mu(p)=x$. Analytically near $p$, the family has equation
$uv=\delta$ for a holomorphic function $\delta$ on $B$.
By the choice of $T_{[D]}B$, the function $\delta$ has a simple
zero at $[D]$. Taking $t=\delta$ as a local parameter on $B$,
we obtain $uv=t$.  After shrinking $B$,
the noncentral fibres are smooth rational curves.
The family $\mathcal{D}_B \subset X \times B$ is a relative local complete intersection and its conormal
bundle  restricts to the
conormal bundle of each fibre in $X$. Semicontinuity and
Lemma~\ref{lem:conormal-vanishing} give
$H^0(C,N_{C/X}^*)=0$ for a nearby fibre $C$. Since
$N_{C/X}$ has rank and degree $2n$, its splitting on $C$
must be $\OO_C(1)^{\oplus 2n}$.

Trivialize $\mu^*L$ near $p$. Since $\theta$ vanishes on the
tangent lines of both central branches, we can write
$\mu^*\theta=v\lambda_1\,du+u\lambda_2\,dv$
for holomorphic functions $\lambda_1,\lambda_2$. This form
vanishes at $p$, and its exterior derivative there gives
$\lambda_2(p)-\lambda_1(p)\ne0$.
In the chosen frame, this difference represents
$-\Levi_x(d\mu_p(\partial_u),d\mu_p(\partial_v))$, which is
nonzero by \eqref{eq:pair}. The vector field
$u\partial_u-v\partial_v$ is tangent to the fibres and
\[
(\mu^*\theta)(u\partial_u-v\partial_v)
=t(\lambda_1-\lambda_2).
\]
Hence $\theta|_{T_C}\ne0$ on a nearby smooth fibre.
As $T_C$ and $L|_C$ both have degree two, $\theta|_{T_C}$ is an
isomorphism. Thus $C$ is a twistor conic.
\end{proof}

\section{The two-point section and homogeneity}\label{sec:section}

\subsection{The two-point evaluation}

Fix a twistor conic $j:C\hookrightarrow X$ and let $\MM$ be the
reduced irreducible component through $[j]$ of the stack of stable maps $\overline{\mathcal M}_{0,0}(X,\beta)$, where $\beta=j_*[C]$.  Denote by
\[
\pi:\UU\to\MM,\qquad e_{\UU}:\UU\to X,\qquad
\II=\UU\times_{\MM}\UU
\]
 its universal curve, evaluation and fibre square. Let
$r_a,r_b:\II\to\UU$ be the projections and let $Y$ be the
coarse moduli space of $\II$. The evaluations
$e_{\UU}\circ r_a$ and $e_{\UU}\circ r_b$ induce
\[
q=(q_a,q_b):Y\longrightarrow X\times X,\qquad
[h:\Sigma\to X;a,b]\longmapsto(h(a),h(b)).
\]
The stack $\MM$ is proper and its coarse moduli space $M_\beta$
is projective.
Let $Y^{\mathrm{tw}}\subset Y$ be the open subset where $h$
embeds $\Sigma$ as a twistor conic.

\begin{lemma}\label{lem:incidence}
(i) $Y$ is an integral projective variety of dimension $4n+2$;
(ii) $Y$ and $q_a$ are smooth along $Y^{\mathrm{tw}}$;
(iii) The morphism $q$ is \'etale at every
$[h;a,b]\in Y^{\mathrm{tw}}$ with $a\ne b$.
It is surjective and generically finite.
\end{lemma}

\begin{proof}
\emph{(i)} Near an embedded twistor conic, $\MM$ agrees with
$\Hilb(X)$, which is smooth of dimension $4n$ at the corresponding
point \cite[Section~3, p.~1065]{IN}. Thus $\dim\MM=4n$.
The morphism $\II=\UU\times_\MM\UU\to\MM$ is flat, with integral
generic fibre $\PP^1\times\PP^1$. Since $\MM$ is integral, so are
$\II$ and its coarse moduli space $Y$, which gives $\dim Y=4n+2$.

An unmarked stable map $h:\Sigma\to X$ of $L$-degree two has no
contracted component. Otherwise such a component would meet the
rest of $\Sigma$ in at least three points by stability, and each of the
resulting branches would contain a noncontracted component,
contrary to the $L$-degree bound. Hence every $h$ has finite fibres.
The proper morphism
\[
Y\longrightarrow M_\beta\times X\times X,\qquad
[h;a,b]\longmapsto([h],h(a),h(b))
\]
therefore has finite fibres and is finite, which shows that $Y$ is projective.

\emph{(ii)} Let $[h:\Sigma\to X;a,b]\in Y^{\mathrm{tw}}$ and put
$C=h(\Sigma)$. Near this point, $Y$ is the fibre square of the
universal family over $\Hilb(X)$, and is therefore smooth.
Since $N_{C/X}\simeq\OO_C(1)^{\oplus 2n}$ is globally generated,
the evaluation $e_H:U_H\to X$ is smooth along $C$.
The projection from the fibre square to $U_H$ is smooth, so
$q_a$ is smooth at $[h;a,b]$.

\emph{(iii)} If $a\ne b$, then
\[
\ker(d_{[h;a,b]}q)
 \simeq H^0\bigl(C,N_{C/X}(-h(a)-h(b))\bigr)=0.
\]
Both $Y$ and $X\times X$ are smooth of dimension $4n+2$ at
these points, so $q$ is \'etale there. In particular, $q$ is
dominant and generically finite. Its image is closed because
$Y$ is projective, hence $q$ is surjective.
\end{proof}

\begin{remark}
For the adjoint variety $X_{\mathfrak g}$, with $\mathfrak g$
neither of type $A$ nor of type $C$, the twistor conics form
the open orbit $G/K\subset\Hilb(X_{\mathfrak g})$, where
$K=\operatorname{Stab}_G(C_0)$ and $C_0$ is the conic in
\eqref{eq:adjoint-conic}
(\cite[Theorem~3.6 and Lemmas~4.2--4.3]{Kwon}).
The universal curve over $G/K$ is $G\times^K C_0$. Hence
\[
Y^{\mathrm{tw}}\simeq
(G\times^K C_0)\times_{G/K}(G\times^K C_0)
\simeq G\times^K(C_0\times C_0),
\]
under which $q|_{Y^{\mathrm{tw}}}$ is given by
$[g;a,b]\mapsto(g\cdot a,g\cdot b)$. Kwon's description of $K$ also shows that $q$ identifies the locus
$a\ne b$ in $Y^{\mathrm{tw}}$ with the open $G$-orbit in
$X_{\mathfrak g}\times X_{\mathfrak g}$, which implies in particular that $q$ is
birational.
\end{remark}

\subsection{Extension of the section}

Let $\MM^\circ\subset\MM$ be the open substack of maps with smooth
domain, and let $\UU^\circ$, $\II^\circ$ and $Y^\circ$ be the
corresponding open subsets. Write $T_\pi$ for the relative tangent
bundle over $\MM^\circ$. The family $\UU^\circ\to\MM^\circ$ is
\'etale locally a product with $\PP^1$. Since $\psi_{\PP^1}$
in \eqref{eq:two-point-tensor} is invariant under automorphisms
of $\PP^1$, it gives a section $\psi_\pi$ of
$r_a^*T_\pi\otimes r_b^*T_\pi$ on $\II^\circ$.
The contraction in \eqref{eq:curvewise-section} gives
\[
\widetilde\eta_{\II}
=\bigl(r_a^*(e_{\UU}^*\theta\circ d_\pi e_{\UU})\otimes
r_b^*(d_\pi e_{\UU})\bigr)(\psi_\pi).
\]
As a morphism to $\Tot(L\boxtimes T_X)$ over $X\times X$, it
descends to a section $\widetilde\eta$ of $q_a^*L\otimes q_b^*T_X$
on $Y^\circ$ by the universal property of the coarse moduli space. We're going to show it extends to $Y$. Let's start with a toy example.

\begin{example}\label{ex:conic-two-point-tensor}
Consider $C_t=\{z_0z_1=t z_2^2\}\subset\PP^2$ and let
$q_t(z)=z_0z_1-tz_2^2$. Its polar form is
\[
Q_t(z_a,z_b)
=q_t(z_a+z_b)-q_t(z_a)-q_t(z_b)
=z_{0,a}z_{1,b}+z_{1,a}z_{0,b}
 -2t z_{2,a}z_{2,b}.
\]
This is a section of
$\OO_{C_t}(1)\boxtimes\OO_{C_t}(1)$. For $t\ne0$ and
fixed $a\in C_t$, the equation $Q_t(a,b)=0$ is the
tangent line to $C_t$ at $a$, since
$Q_t(z_a,z_b)=dq_t|_{z_a}(z_b)$. That line meets $C_t$
in the divisor $2a$. Hence $Q_t$ vanishes precisely
twice along the diagonal of $C_t\times C_t$.

For $t\ne0$, the adjunction isomorphism
$\iota_t:\OO_{C_t}(1)\simeq T_{C_t}$ sends $Q_t$ to a
multiple of $\psi_{C_t}$. To determine the multiple, we work on
$z_2\ne0$ and put $u=z_0/z_2$ and $v=z_1/z_2=t/u$.
Normalize $\iota_t$ so that
$\iota_t(z_2)=u\partial_u-v\partial_v$, which acts as
$u\partial_u$ in the coordinate $u$. Then
\[
\frac{Q_t}{z_{2,a}\otimes z_{2,b}}
=u_av_b+v_au_b-2t
=(u_b-u_a)(v_a-v_b)
=\frac{t(u_b-u_a)^2}{u_au_b},
\]
and hence
\[
(\iota_{t,a}\otimes\iota_{t,b})(Q_t) =\frac{t(u_b-u_a)^2}{u_au_b}
  (u_a\partial_{u_a})\otimes(u_b\partial_{u_b})
  =t(u_b-u_a)^2 \partial_{u_a}\otimes\partial_{u_b}
=t\psi_{C_t}.
\]
This identity extends from the chart to $C_t\times C_t$.
Although $\psi_{C_t}$ is defined here only for $t\ne0$,
$Q_t$ remains regular at $t=0$. This regular tensor is used
below to extend $\widetilde\eta$ across the boundary.
\end{example}

\begin{proposition}\label{prop:global-section}
The section $\widetilde\eta$ on $Y^\circ$ extends uniquely to a
section of $q_a^*L\otimes q_b^*T_X$ on $Y$.
\end{proposition}

\begin{proof}
The section is already defined over the smooth-domain locus.
Since regularity can be checked \'etale locally, we first extend
it over a reduced affine \'etale chart $S=\Spec R$ of $\MM$ near
a nodal map, then compare the local extensions on overlaps.

Let $\pi_S:\mathcal U_S\to S$ and $e_S:\mathcal U_S\to X$
be the universal curve and evaluation. Put $\mathcal L_S=e_S^*L$.
By the proof of Lemma \ref{lem:incidence}, each geometric fibre $C_s$ is either a smooth $\PP^1$ with
$\mathcal L_s\simeq\OO_{\PP^1}(2)$, or two copies of $\PP^1$
meeting at one node, with $\mathcal L_s$ of multidegree $(1,1)$.
The normalization sequence on a nodal fibre, and the corresponding
calculation on a smooth fibre, give
\[
H^1(C_s,\mathcal L_s^m)=0,\qquad
h^0(C_s,\mathcal L_s^m)=2m+1
\quad (m=1,2).
\]
In particular, $\mathcal L_s$ is very ample on every fibre.

By cohomology and base change,
$\mathcal{E}:=\pi_{S*}\mathcal L_S$ and
$\mathcal{F}:=\pi_{S*}\mathcal L_S^2$ are locally free of ranks three and
five. After shrinking $S$ and trivializing $\mathcal{E}$, evaluation gives
a morphism $j:\mathcal U_S\to\PP^2_R$. On each geometric fibre
it is the complete conic embedding: a smooth fibre maps to a
smooth conic, and the two components of a nodal fibre map to
distinct lines. The morphism $j$ is a closed immersion. Indeed,
it is proper with finite fibres and hence finite; fibrewise
surjectivity of
$\OO_{\PP^2_R}\to j_*\OO_{\mathcal U_S}$ then gives surjectivity
by Nakayama's lemma.

The multiplication map
$\operatorname{Sym}^2 \mathcal{E}\longrightarrow \mathcal{F}$
is surjective on every fibre, with kernel generated by the
equation of its conic. Thus its kernel is a line bundle. After
shrinking $S$ once more, a generator gives a quadratic form
$q\in R[z_0,z_1,z_2]_2$. Its zero scheme is flat over $S$ and
agrees with $j(\mathcal U_S)$ on every geometric fibre, so the
two schemes agree. In particular,
$\mathcal L_S=\OO_{\mathcal U_S}(1)$.

Near a nodal fibre, choose a smooth point and, after an \'etale
base change, a section of $\pi_S$ through it. Choose coordinates
placing this section at $[1:0:0]$ and its relative tangent line
at $z_1=0$. After multiplying $q$ by a unit, it has the form
\[
q=z_0z_1+\lambda z_1^2+\mu z_1z_2+\nu z_2^2.
\]
Replacing $z_0$ by $z_0+\lambda z_1+\mu z_2$ and putting
$\delta=-\nu$ gives
\[
\mathcal U_S=\{z_0z_1=\delta z_2^2\}\subset\PP^2_R,
\qquad
\mathcal L_S=\OO_{\mathcal U_S}(1).
\]

Let $D(\delta):=\{s\in S| \delta(s)\ne0\}=S\setminus V(\delta)$.
By the equation $z_0z_1=\delta z_2^2$, this is the
smooth-domain locus of $\pi_S$.  It meets every irreducible component of the reduced scheme $S$.
 Hence $\delta$ is a nonzerodivisor
in $R$ and remains a nonzerodivisor on the flat fibre square
$\mathcal U_S\times_S\mathcal U_S$.

Adjunction identifies $\omega_{\mathcal U_S/S}^{-1}$ with
$\mathcal L_S$. Dualizing the natural map
$\Omega^1_{\mathcal U_S/S}\to\omega_{\mathcal U_S/S}$
gives $\iota:\mathcal L_S\longrightarrow T_{\mathcal U_S/S}$.
On $z_2\ne0$, set $u=z_0/z_2$ and $v=z_1/z_2$, so that
$uv=\delta$. We normalize the identification by
$\iota(z_2)=u\partial_u-v\partial_v.$
Since $\pi_{S*}\OO_{\mathcal U_S}=\OO_S$, 
we have
\[
\Hom_{\mathcal U_S}(\mathcal L_S,\mathcal L_S)
=H^0(\mathcal U_S,\OO_{\mathcal U_S})
=H^0(S,\OO_S)=R.
\]
Then the composite
\[
\mathcal L_S\xrightarrow{\iota}T_{\mathcal U_S/S}
\xrightarrow{d_{\pi_S}e_S}e_S^*T_X
\xrightarrow{e_S^*\theta}\mathcal L_S
\]
is multiplication by a function $\gamma\in R$.
On the node subscheme $V(u,v)\simeq\Spec(R/(\delta))$,
both coefficients of $\iota(z_2)$ vanish, so its image under
$d_{\pi_S}e_S$ restricts to zero there. Since $z_2$ trivializes
$\mathcal L_S$ on this subscheme, we have $\gamma\in\delta R$.

Let $Q$ be the section of
$\mathcal L_{S,a}\otimes\mathcal L_{S,b}$ obtained from
Example~\ref{ex:conic-two-point-tensor} by replacing $t$
with $\delta$. Thus
\[
Q=z_{0,a}\otimes z_{1,b}+z_{1,a}\otimes z_{0,b}
  -2\delta z_{2,a}\otimes z_{2,b}.
\]
The calculation in that example gives
\[
(\iota_a\otimes\iota_b)(Q)=\delta\psi_\pi
\quad\text{over }D(\delta).
\]
Using the definition of $\widetilde\eta_{\II}$ and the identity
$e_S^*\theta\circ d_{\pi_S}e_S\circ\iota=\gamma\,\mathrm{id}$,
we obtain on $D(\delta)$
\[
\widetilde\eta_{\II}
=\frac{\gamma}{\delta}
\bigl(\mathrm{id}\otimes(d_{\pi_S}e_S\circ\iota)_b\bigr)(Q).
\]
Because $\gamma/\delta\in R$, the right-hand side is regular
on the whole fibre square.

On overlaps, the local extensions agree over the smooth-domain
locus. They agree everywhere because $\delta$ is a
nonzerodivisor on the fibre square and the target is locally
free. They therefore glue to a section on $\II$.
Viewed as a morphism
$\II\to\Tot(L\boxtimes T_X)$ over $X\times X$, this section
factors through the coarse moduli space $Y$, giving the required
extension. Its uniqueness follows from the integrality of $Y$
and the density of $Y^\circ$.
\end{proof}

The relative anticanonical embedding of smooth or one-nodal
rational curves as plane conics is proved in
\cite[Lemma~8 and the proof of Proposition~6]{EF}.
In the precedent proof, we obtain the
plane-conic model directly from $e_S^*L$ and the pole cancellation  in this model parallels \cite[Proposition~3.3]{HLW}.

\subsection{Descent to $X\times X$  and the quadratic term}

Let $\mathcal I_x$ be the ideal sheaf of $x \in X$ and
$\mathfrak m_x$ its stalk at $x$. For a section $s$
of a vector bundle $E$ vanishing to order at least two at $x$, write
\[
[s]_{2,x}\in(\mathfrak m_x^2/\mathfrak m_x^3)\otimes E_x
\simeq\operatorname{Sym}^2T_x^*X\otimes E_x
\]
for its quadratic term, evaluated as a homogeneous quadratic
polynomial. This can be computed by taking a curve $c(u)$ with $c(0)=x$ and $c'(0)=v$, then the Taylor expansion gives
$s(c(u))=u^2[s]_{2,x}(v)+O(u^3)$ (under a local
trivialization of $E$).

\begin{proposition}\label{prop:section}
Let $x\in X$ and let $Y_x$ be an irreducible component of
$q_a^{-1}(x)_{\mathrm{red}}$ meeting $Y^{\mathrm{tw}}$.
Put $q_x=q_b|_{Y_x}$ and use the canonical identification
$q_a^*L|_{Y_x}=L_x\otimes\OO_{Y_x}$. Then we have:

(i) $Y_x$ is integral and projective, and $q_x$ is
surjective and generically finite;

(ii) There is a unique section
$\xi_x\in L_x\otimes H^0(X,T_X)$ such that
$q_x^*\xi_x=\widetilde\eta|_{Y_x}$;

(iii) The section $\xi_x$ vanishes to order at least two
at $x$, and for $v \in T_xX$, we have 
$[\xi_x]_{2,x}(v)=\theta_x(v)\otimes v$.
\end{proposition}

\begin{proof}
\emph{(i)} As a component of a reduced closed subscheme of $Y$,
the variety $Y_x$ is integral and projective. Choose
$[h;a,b]\in Y_x\cap Y^{\mathrm{tw}}$ and put $C=h(\Sigma)$.
Let $B_x$ be an irreducible open neighbourhood of $[C]$ in
$H_x\cap H^{\mathrm{tw}}$, small enough for the identification
between embedded maps in $\MM$ and their images in $\Hilb(X)$.
Then
\[
U_H\times_H B_x\longrightarrow q_a^{-1}(x),\qquad
([C'],y)\longmapsto[C'\hookrightarrow X;x,y]
\]
identifies its source with an irreducible open subset containing
$[h;a,b]$. By Lemma~\ref{lem:incidence}, the fibre
$q_a^{-1}(x)$ is smooth along this subset. It therefore lies
in the component $Y_x$.

At points with $y\ne x$, the map $q_x$ is \'etale: locally on
this subset it is the base change of the \'etale map $q$ along
$\{x\}\times X\hookrightarrow X\times X$.
Thus $q_x$ has a nonempty open image and
$\dim Y_x=\dim X$. It is dominant and generically finite.
Its image is closed by projectivity, so it is surjective.

\emph{(ii)} Proposition~\ref{prop:global-section} gives
$\widetilde\eta|_{Y_x}\in
H^0(Y_x,L_x\otimes q_x^*T_X)$.
Since $q_x$ satisfies the hypotheses of
Lemma~\ref{lem:descent}, tensoring its pullback isomorphism
with $L_x$ gives the unique section $\xi_x$.

\emph{(iii)} Parametrize a curve in $B_x$ by
$h':\PP^1\hookrightarrow X$, and choose an affine coordinate
$u$ with $h'(0)=x$. Put $v=dh'_0(\partial_u)$.
Formula~\eqref{eq:curvewise-section} gives
\[
\xi_x(h'(u))
=u^2\theta_x(v)\otimes dh'_u(\partial_u).
\]
Hence $\xi_x(x)=0$ and $d_x\xi_x(v)=0$.
By Lemma~\ref{lem:twistor-tangent-map}, the directions $[v]$
obtained from $B_x$ contain a nonempty open subset of
$\PP(T_xX)$. Since $d_x\xi_x$ is linear, it vanishes on
$T_xX$.

In a local trivialization near $x$,
$dh'_u(\partial_u)=v+O(u)$. The coefficient of $u^2$ in the
precedent formula is therefore $\theta_x(v)\otimes v$.
This proves $[\xi_x]_{2,x}(v)=\theta_x(v)\otimes v$ for the directions
obtained from $B_x$. Both sides are homogeneous quadratic
maps on $T_xX$, so the equality holds for every $v\in T_xX$.
\end{proof}

\begin{lemma}\label{lem:birational-evaluation}
The morphism $q:Y\to X\times X$ is birational.
\end{lemma}

\begin{proof}
Let $\Omega=\{[h;a,b]\in Y^{\mathrm{tw}}:a\ne b\}.$
Since $Y$ is integral and $\dim Y=\dim(X\times X)$,
the closed subset $Y\setminus\Omega$ has image of smaller
dimension. As $q$ is proper, $q(Y\setminus\Omega)$ is a proper
closed subset of $X\times X$. Choose a general $(x,y)$ outside
this image. The fibre $q^{-1}(x,y) \subset \Omega$ is nonempty and finite.

Choose $z\in q^{-1}(x,y)$. Since $q_a$ is smooth at $z$,
the point $z$ lies on a unique irreducible component $Y_x$
of $q_a^{-1}(x)_{\mathrm{red}}$. Proposition~\ref{prop:section}
gives $\xi_x\in L_x\otimes H^0(X,\mathcal I_x^2\otimes T_X)$.
By \eqref{eq:curvewise-section},
$\xi_x(y)=\widetilde\eta(z)\ne0$. Contracting $\xi_x$ with
a nonzero element of $L_x^*$ gives a vector field $\mathbf{v}$ with
a double zero at $x$ and $\mathbf{v}(y)\ne0$.

For any twistor conic $D$ through $x$ and $y$,  we have
$\mathbf{v}|_D\in H^0(D,T_X|_D(-2x))$ since $\mathbf{v}$ has double zero
at $x$. Its image under
$T_X|_D\to N_{D/X}$ is zero, since
\[
H^0(D,N_{D/X}(-2x))
=H^0\bigl(\PP^1,\OO_{\PP^1}(-1)^{\oplus 2n}\bigr)=0.
\]
Thus the vector field $\mathbf{v}$ is tangent to $D$.
Since $\mathbf{v}(y)\ne0$, every such $D$
contains the unique germ of an integral curve of $\mathbf{v}$ through
$y$. The conics therefore have the same germ at $y$ and
coincide. It follows that  $q^{-1}(x,y)$
has one point. The fibre is reduced because $q$ is \'etale
on $\Omega$. Hence  $q$ is birational.
\end{proof}

\begin{proposition}\label{prop:diagonal-section}
(i) There is a unique section
$\eta\in H^0(X\times X,L\boxtimes T_X)$ such that
$q^*\eta=\widetilde\eta$.

(ii) For every $x$ and every $Y_x$ as in
Proposition~\ref{prop:section}, one has $\eta(x,-)=\xi_x$.
In particular, $\xi_x$ is independent of the choice of $Y_x$.

(iii)  Let $\Delta\subset X\times X$ be the diagonal, with ideal sheaf
$\mathcal I_\Delta$. We identify $N_{\Delta/(X \times X)}$ with $T_X$ by 
$(v_1, v_2) \mapsto v_2-v_1$.  Then the section $\eta$ satisfies
$\eta\in H^0\left(X\times X,
\mathcal I_\Delta^2\otimes(L\boxtimes T_X)\right)$.  Write $[\eta]_{2,x}$ for the value at $(x,x)$ of its quadratic class modulo $\mathcal I_\Delta^3$. Then for  $x\in X$ and  $v\in T_xX$ we have
$$[\eta(-,x)]_{2,x}(v)=[\eta]_{2,x}(-v)
=[\eta(x,-)]_{2,x}(-v)
=\theta_x(v)\otimes v.$$
\end{proposition}

\begin{proof}
(i) By Lemma~\ref{lem:birational-evaluation},
$q$ is birational. It is proper, $Y$ is integral, and
$X\times X$ is normal; hence $q_*\OO_Y=\OO_{X\times X}$.
The projection formula gives
$q_*q^*(L\boxtimes T_X)\simeq L\boxtimes T_X.$
Taking global sections yields the unique section $\eta$.

(ii) Restriction to $Y_x$ gives
$q_x^*\eta(x,-)=\widetilde\eta|_{Y_x}=q_x^*\xi_x.$
Since $q_x$ is surjective, $\eta(x,-)=\xi_x$.

(iii)
Put $U=q_a(Y^{\mathrm{tw}})$. By Lemma~\ref{lem:incidence},
$U$ is a nonempty open subset of $X$. For $x\in U$,
part~\textup{(ii)} and Proposition~\ref{prop:section} show that
$\eta(x,-)$ vanishes to order at least two at $x$.
Thus $\eta|_\Delta$ vanishes on $U$, and hence on all of $\Delta$.

Choose local coordinates $z_i$ on $X$ and put
$t_i=z_i(y)-z_i(x)$. The class of $\eta$ modulo
$\mathcal I_\Delta^2$ is its part linear in the $t_i$;
at $(x,x)$ this is the linear term of $\eta(x,y)$ at $y=x$.
It vanishes for $x\in U$, hence everywhere. Therefore
$\eta\in H^0(X\times X,\mathcal I_\Delta^2\otimes (L\boxtimes T_X))$.

The class of $\eta$ modulo $\mathcal I_\Delta^3$ is now a
section of
\[
(\mathcal I_\Delta^2/\mathcal I_\Delta^3)\otimes (L\boxtimes T_X)|_\Delta
\simeq \operatorname{Sym}^2\Omega_X^1\otimes L\otimes T_X.
\]

The inclusions $X\times\{x\}$ and $\{x\}\times X$
send $v\in T_xX$ to $-v$ and $v$, respectively, in
$N_{\Delta/(X\times X)}$. This gives the first two
equalities in (iii).
For $x\in U$, part~\textup{(ii)}
and Proposition~\ref{prop:section} give
$[\eta]_{2,x}(v)=\theta_x(v)\otimes v.$
Both sides define sections of the same bundle, so their
equality on $U$ extends to all of $X$.
\end{proof}

The curvewise formula for $\widetilde\eta$ agrees with that of $W_{21}$ in
\cite[Proposition~5.1]{OpenAI}.  In
\cite[Proposition~6.1 and Corollary~6.2]{OpenAI}, it proves that the contact
contraction of $W_{21}$ descends as a symmetric section of $L\boxtimes L$
and it uses its cubic term to prove the homogeneity
(\cite[Lemma~7.3 and Theorem~7.5]{OpenAI}).

\begin{proof}[Proof of Theorem~\ref{thm:main}]
%If $\dim X=1$, then $(X,L)\simeq(\PP^1,\OO(2))$, the adjoint variety of $\mathfrak{sl}_2$. 
By the reduction in
Section~\ref{sec:lines}, we may  assume
$\Pic(X)=\mathbb Z[L]$.

Let $\eta$ be the section of
Proposition~\ref{prop:diagonal-section}. For $x\in X$, put
\[
E_x=\operatorname{im}\bigl(H^0(X,T_X)\longrightarrow T_xX\bigr).
\]
By the K\"unneth formula, $\eta$ is a sum of sections
$s\boxtimes w$, with $s\in H^0(X,L)$ and $w\in H^0(X,T_X)$.
Thus $\eta(-,x)$ takes values in $L\otimes E_x$ and its
quadratic term at $x$ takes values in $L_x\otimes E_x$.
By Proposition \ref{prop:diagonal-section}(iii), we have
$[\eta(-,x)]_{2,x}(v) =\theta_x(v)\otimes v
\in L_x\otimes E_x.$
 It follows that $v\in E_x$ whenever $\theta_x(v)\ne0$. Such vectors  span $T_xX$, hence  $E_x=T_xX$ for every $x$, so $T_X$ is globally generated.

It follows that
$\operatorname{Aut}^0(X)$ acts transitively on $X$. Since every
global vector field preserves $F$,
$X$ is a homogeneous complex contact manifold. As $X$ is
Fano, it is simply connected, so Boothby's classification
\cite[\S2, pp.~1034--1036]{Wolf} gives a contact isomorphism
$X\simeq X_{\mathfrak g}$ identifying the contact line bundles.

\end{proof}

 State Key Laboratory of Mathematical Sciences, Morningside Center of Mathematics, Academy of Mathematics and Systems Science, Chinese Academy of Sciences, Beijing 100190, China; and School of Mathematical Sciences, University of Chinese Academy of Sciences, Beijing, China \quad  (bhfu@math.ac.cn)
\end{document}